\documentclass[reqno, twoside]{amsart}

\newtheorem{theorem}{Theorem}[section]
\newtheorem{lemma}[theorem]{Lemma}

\theoremstyle{definition}
\newtheorem{definition}[theorem]{Definition}

\theoremstyle{remark}

\numberwithin{equation}{section}

\usepackage[english]{babel}
\usepackage{amsfonts, amsmath}
\usepackage{graphics}
  \usepackage{amssymb}
\usepackage{amsmath}
  \usepackage{latexsym}
  \usepackage{amsthm}
\usepackage[all]{xy}

\def\H{\mathord{\mathbb H}}
\def\C{\mathord{\mathbb C}}
\def\R{\mathord{\mathbb R}}

\def\3{{\ss}}
\def\2{\frac{1}{2}}
\def\4{\frac{1}{4}}
\def\8{\frac{1}{8}}
\def\x{\times}
\def\.{\cdot}

\def\<{\langle}
\def\>{\rangle}

\def\Im{\mathop{\rm Im\,}\nolimits}
\def\Re{\mathop{\rm Re\,}\nolimits}

\def\Ad{\mathop{\rm Ad}\nolimits}

\def\trace{\mathop{\rm trace\,}\nolimits}

\def\g{\mathfrak{g}}

\def\p{\mathfrak{p}}

\def\k{\mathfrak{k}}
\def\q{\mathfrak{q}}

\def\h{\mathfrak{h}}

\def\ms{\medskip\noindent}
\def\bsm{\begin{smallmatrix}}
\def\esm{\end{smallmatrix}}
\def\bpm{\begin{pmatrix}}
\def\epm{\end{pmatrix}}
\def\beq{\begin{equation}}
\def\eeq{\end{equation}}

\begin{document}

\title{Almost positive curvature\\ 
on the Gromoll-Meyer sphere}

\author{J.-H. Eschenburg}
\address{Institut f\"ur Mathematik, Universit\"at Augsburg,
D-86135 Augsburg, Germany}
\email[(Eschenburg)]{eschenburg@math.uni-augsburg.de}

\author{M. Kerin}
\address{Department of Mathematics, University of Pennsylvania, 209 S 33rd St., Philadelphia, PA 19104, USA}
\email[(M. Kerin)]{mkerin@math.upenn.edu}

\date{\today}

\thanks{The second author would like to thank the University of Pennsylvania for financial support.}



\begin{abstract} Gromoll and Meyer have represented a certain exotic 7-sphere
$\Sigma^7$ as a biquotient of the Lie group $G = Sp(2)$. We show for a
2-parameter family of left invariant metrics on $G$ that the induced metric
on $\Sigma^7$ has strictly positive sectional curvature at all points outside
four subvarieties of codimension $\geq 1$ which we describe explicitly.
\end{abstract}

\maketitle

\section{Introduction}

Let $G = Sp(2)$ be the Lie group of unitary quaternionic $2\x2$-matrices.
Consider the subgroup $U \subset G\x G$,  \beq \label{U}  U =
\{\left((\bsm q & \cr & 1 \esm),(\bsm q & \cr & q\esm)\right); \ q \in
Sp(1)\},  \eeq  which acts on $G$ by left and right translations.  D.
Gromoll and W. Meyer \cite{GM} have shown that this action is free and
that the orbit space $M = G/U$ is a smooth manifold which is an exotic
7-sphere (homeomorphic but not diffeomorphic to the standard 7-sphere).
If $G$ is equipped with a Riemannian metric of nonnegative  sectional
curvature whose isometry group contains $U$, then  by O'Neill's formula
\cite{B} the orbit space $M = G/U$ inherits a Riemannian  metric of
nonnegative sectional curvature. Thus starting with the bi-invariant
metric on $G$, Gromoll and Meyer constructed a metric of nonnegative
sectional curvature on the exotic sphere $M$. In fact the curvature is
strictly positive on some nonempty open subset of $M$. However, as was
observed by F. Wilhelm \cite{W}, there is also an open subset with zero
curvature planes in the tangent space of each of its points. But Wilhelm
constructed another $U$-invariant metric on $Sp(2)$ (neither left nor
right invariant) for which the curvature of $M$ is strictly positive
outside a subset of measure zero in $M$ (``almost positive curvature'').
In \cite{E2} the same fact was claimed for a much simpler and left
invariant metric on $Sp(2)$; however, as was pointed out by the  second
author, the proof contains a serious mistake (see Remark 3 at the end of
the present paper). The  purpose of our paper is to correct this error.
In fact we prove  the following result, some ideas of which go back to
\cite{E1} (see  Theorem \ref{theorem} for details):

\begin{theorem} There is a left invariant and $U$-invariant metric on  $G
= Sp(2)$ such that the induced metric on $M = G/U$ has strictly  positive
curvature outside a finite union of subvarieties of codimension $\geq 1$.
The zero curvature set $Z \subset M$ can be explicitly determined.
\end{theorem}

\section{Cheeger metrics on Lie groups}

On each Riemannian manifold, let us denote
\begin{eqnarray}
	\kappa(X,Y) &=& \<R(X,Y)Y,X\>,	\cr
	\sec(X,Y) &=& \kappa(X,Y)/|X\wedge Y|^2
\end{eqnarray}
for any two tangent vectors $X,Y$; the second expression is the 
sectional curvature of the plane $\sigma$ spanned by $X,Y$.

Let $G$ be a Lie group with a left invariant metric $\<\ ,\ \>$ of nonnegative
sectional curvature. Suppose that the metric is also right invariant with
respect to a compact subgroup $K \subset G$, hence the induced metric on
$K$ is bi-invariant. The Lie algebras of $G$ and $K$ will be denoted $\g$ and
$\k$. We may contract the metric on $G$ in the direction
of the $K$-cosets by viewing $G$ as the homogeneous space $(G\x K)/\Delta
K$ (where $\Delta K = \{(k,k);\ k\in K\}$) and choosing the metric
induced from the Riemannian product metric on $G\x sK$ ({\it Cheeger
contraction}, cf. \cite{C}, \cite{B}) where $sK$ is $K$ with  metric scaled
by $s > 0$. A vector $(X,X') \in \g \x \k$ is perpendicular to the 
$\Delta K$-orbit (``horizontal'') iff $X+sX' \perp \k$, i.e. $X' =
- s^{-1}X_\k$ where $X_\k$ is the  $\k$-projection of $X$. Using the
Riemannian submersion $G\x K \to G$, $(g,k) \mapsto gk^{-1}$, a
horizontal vector $(X,-s^{-1}X_\k)\in\g\x\k$ is mapped onto 
$X+s^{-1}X_\k = X_\perp + (1+s^{-1})X_\k \in \g$ where $X_\perp = X-X_\k
\in \k^\perp$. Vice versa, the horizontal lift of $X = X_\perp + X_\k \in
\g$ is the horizontal vector
\begin{eqnarray} \label{tildeX}
	\widehat X &=& (\tilde X,-s^{-1}\tilde X_\k),\ \ \textrm{ where}\cr
	\tilde X &=& X_\perp + \frac{s}{s+1}X_\k.
\end{eqnarray}
Thus the new (left invariant) metric is
\begin{eqnarray} \label{(XY1)}
	\<X,Y\>_1 
	&=& \<\widehat X,\widehat Y\> \cr
	&=& \<\tilde X,\tilde Y\> + s\,\<s^{-1}\tilde X_\k,s^{-1}\tilde Y_\k\>  \cr 
	&=& \<\tilde X,\tilde Y\> + s^{-1}\<\tilde X_\k,\tilde Y_\k\>  \cr 
	&=& \<\tilde X_\perp,\tilde Y_\perp\> 
	+ s^{-1}(s+1)\<\tilde X_\k,\tilde Y_\k\> 	\cr
	&=& \<X_\perp,Y_\perp\> + s(s+1)^{-1}\<X_\k,Y_\k\> \cr
	&=& \<\tilde X,Y\>.
\end{eqnarray}
For the curvature terms we have
\beq
	\kappa(\widehat X,\widehat Y) 
	= \kappa(\tilde X,\tilde Y) + s^{-3}\kappa(\tilde X_\k,\tilde Y_\k).
\eeq
Since all terms are nonnegative, the left hand side vanishes if and
only if both summands on the right are zero. Thus a plane $\sigma$ 
spanned by $X,Y \in \g$ has zero curvature in the new metric,
$\sec_1(\sigma) = 0$, if and
only if  $\sec(\tilde\sigma) = 0$ and $[X_\k,Y_\k] = 0$.\footnote
	{The ``if'' statement is not obvious because of the
	nonnegative O'Neill term. However, in all our examples starting
	with a bi-invariant metric on some Lie group, the vanishing of the
	curvature implies that the O'Neill term also vanishes, see
	\cite{E1}, p. 29f, Equations (1) - (4) or \cite{Wi}, \cite{T}}

\ms {\bf Example 1.} Suppose that the initial metric $\<\ ,\ \>$ 
on $G$ is bi-invariant. Let $\g = \k + \p$ be the orthogonal decomposition. 
Consider the above metric 
\beq \label{s}
	\<X,Y\>_1 = \<X_\p,Y_\p\> + \tilde s\<X_\k,Y_\k\>
\eeq
with $\tilde s = \frac{s}{s+1}$. 
Then $\sec(\tilde\sigma) = 0$ $\iff$ $[\tilde X,\tilde Y] =
0$, and hence $\sec_1(\sigma) = 0$ $\iff$ 
$$
	[\tilde X,\tilde Y] = 0,\ \ \ [X_\k,Y_\k] = 0. 
$$
If $(G,K)$ is a symmetric pair, i.e. the orthogonal
complement $\p \subset\g$ satisfies $[\p,\p] \subset \k$, then
$|\tilde X,\tilde Y]_\k = [\tilde X_\k,\tilde Y_\k] + [\tilde X_\p,\tilde Y_\p]$ 
and $|\tilde X,\tilde Y]_\p = [\tilde X_\k,\tilde Y_\p] + [\tilde X_\p,\tilde
Y_\k]$,  hence $\sec_1(\tilde\sigma) = 0$ $\iff$ 
\beq 
	0 = [X_\k,Y_\k] = [X_\p,Y_\p] = [X_\k,Y_\p] + [X_\p,Y_\k] = [X,Y].
\eeq

\ms  {\bf Example 2.} Let $G \supset K \supset H$ a chain of subgroups
and suppose that both $(G,K)$ and $(K,H)$ are symmetric pairs. Let $\g =
\k + \p$ and $\k = \h + \q$ be the corresponding decompositions.  Start
with the metric $\<\ ,\ \>_1$ defined by Example 1, depending
on a parameter $s > 0$, and define the metric $\<\ ,\ \>_2$ by Cheeger
contraction along $H$ (depending on a new parameter $t>0$) as in (\ref{(XY1)})
where $K$ is replaced by $H$ and $\<\ ,\ \>_1$ takes the role of $\<\ ,\ \>$:
\goodbreak\begin{eqnarray} \label{t}
  \<X,Y\>_2 &=& \<X_\p,Y_\p\>_1 + \<X_\q,Y_\q\>_1 + \tilde t\<X_\h,Y_\h\>_1 \cr
&=& \<X_\p,Y_\p\> + \tilde s\<X_\q,Y_\q\> 
+ \tilde s\tilde t\<X_\h,Y_\h\>
\end{eqnarray}
with $\tilde t = \frac{t}{t+1}$. Then $\sec_2(\sigma) = 0$ $\iff$
$\sec_1(\tilde\sigma) = 0$ and $[\tilde X_\h,\tilde Y_\h] = 0$ $\iff$
\beq \label{tildeXY}
0 = [\tilde X,\tilde Y] = [\tilde X_\k,\tilde Y_\k] = [X_\p,Y_\p] 
= [X_\q,Y_\q] = [X_\h,Y_\h], 
\eeq
where $\tilde X = X_\p + X_\q + \frac{t}{t+1}X_\h$ and $\tilde Y = 
Y_\p + Y_\q + \frac{t}{t+1}Y_\h$ like in (\ref{tildeX}).

\section{Zero curvature planes on $Sp(2)$}

Let us consider the chain $G \supset K \supset H$ for  $G = Sp(2)$, $K =
Sp(1)\x Sp(1)$ and $H = \Delta Sp(1) = \{\left(\bsm q & \cr & q \esm
\right);\ q\in Sp(1)\}$. The pairs $(G,K)$ and $(K,H)$ are symmetric,
corresponding to the rank-one symmetric spaces $S^4$ and $S^3$. We start
with the bi-invariant trace metric $\<X,Y\> = \Re\trace X^*Y = \Re\sum
\overline{x_{ij}}\,y_{ij}$ on $\g = \mathfrak{sp}(2)$, apply Cheeger
contraction in the $K$-direction and Cheeger-contract again in the $H$-direction,
defining metrics $\<\ ,\ \>_1$ and $\<\ ,\ \>_2$  as in Example 2. 

Since $G/K = S^4$ as well as $K/H = S^3$ and 
$H = S^3$ have positive curvature,  
the vanishing of the last three brackets in (\ref{tildeXY}) means the
linear dependence of the factors. In particular we may assume 
$Y_\p = 0$, i.e. 
$
	\tilde Y = \tilde Y_\k = \bpm y_1 & \cr & y_2 \epm.
$

\ms {\bf Case 1.} $X_\p = 0$, i.e. 
$\tilde X = \tilde X_\k = \bpm x_1 & \cr  & x_2 \epm$. 

\ms From $[\tilde X_\k,\tilde Y_\k] = 0$ we obtain that the imaginary
quaternions $x_1,y_1$ as well as $x_2,y_2$ are linearly dependent.
Moreover, from $[X_\q,Y_\q] = [X_\h,Y_\h] = 0$ we see that also $x_1\pm
x_2$ and $y_1\pm y_2$ are linearly dependent. Putting $y = y_1$, we may assume
\beq \label{XY1}
	\tilde Y =  \bpm y &  \cr  & 0 \epm,\ \ \ 
	\tilde X =  \bpm 0 &  \cr  & y \epm.
\eeq

\ms
{\bf Case 2.} $X_\p \neq 0$, i.e. 
$X = \bpm x_1 & -\bar x \cr x & x_2 \epm$ with $x \neq 0$:

\ms
Then $0 = [\tilde X,\tilde Y]_\p = [X_\p,\tilde Y]$ $\iff$ $y_2 =
xyx^{-1}$ for $y := y_1$,\\
and $\ \, 0 = [\tilde X,\tilde Y]_\k = [\tilde X_\k,\tilde Y_\k]$ $\iff$
$x_1 = \alpha y_1$, $x_2 = \beta y_2$ for real numbers $\alpha,\beta$, hence
\beq \label{XY0}
	\tilde Y = \bpm  y &  \cr  & xyx^{-1} \epm,\ \ \ 
	\tilde X = \bpm \alpha y & -\bar x \cr x & -\alpha xyx^{-1} \epm
\eeq
where $x,y \in \H$, $y$ imaginary and $\alpha\in \R$;
we have $\beta = -\alpha$ since we require $\<\tilde X,\tilde Y\> = 0$.

\ms {\bf Case 2a.} $\alpha = 0$, hence
\beq \label{XY2}
	\tilde Y =  \bpm y &  \cr  & xyx^{-1} \epm,\ \ \ 
	\tilde X =  \bpm  & -\bar x \cr x &  \epm.
\eeq

\ms {\bf Case 2b.} $\alpha \neq 0$, hence (without loss of generality) 
$\alpha = 1$.

\ms Then $[X_\h,Y_\h] = 0$ iff $y+xyx^{-1}$ and $y-xyx^{-1}$ are
proportional which means $xyx^{-1} = \beta y$.
Comparing the norms on both sides we get 
\beq
	xyx^{-1} = \pm\,y,
\eeq
and
\beq \label{XY3}
	\tilde Y = Y_\pm = \bpm y &  \cr  & \pm y \epm,\ \ \ 
	\tilde X = X_\pm = \bpm y & -\bar x \cr x & \mp y \epm.
\eeq

\begin{lemma} \label{lemma0} The zero curvature planes in $\g = T_eG$ for $G =
Sp(2)$ and the metric $\<\ ,\ \>_2$  are spanned by $X,Y \in \g$ with
$\tilde X,\tilde Y$ given by either (\ref{XY1}) or (\ref{XY2}) or
(\ref{XY3}). \end{lemma}

\section{The Gromoll-Meyer sphere} 

The metric $\<\ ,\ \>_2$ on $G = Sp(2)$ is invariant under the action  of
$U$ (cf. (\ref{U})) and hence it induces a metric on the orbit space $M =
G/U$. Consider any 
\beq
	g = \bpm a & b \cr c & d \epm \in G.
\eeq
Since $g$ is unitary, the rows and columns are unit vectors, in
particular
\beq \label{a2+b2}
	|a|^2+|b|^2 = 1.
\eeq
The vertical space at $g$ of the submersion $\pi : G \to G/U$
is $T_g(U.g) = gV_g$ with $V_g = \{v_g;\ v \in \Im\H\}$ where
\beq
	v_g = g^{-1}\, \bpm v & 0 \cr 0 & 0 \epm\, g - 
	\bpm v & 0 \cr 0 & v \epm =
	\bpm \bar ava-v & \bar avb \cr \bar bva & \bar bvb-v\epm
\eeq
Thus according to (\ref{(XY1)}), a vector $gX \in T_gG$ is horizontal for
$\pi$ iff 
\beq
	0 = \<X,v_g\>_2 = \<\tilde X,v_g\>_1
\eeq
for all $v \in \Im\H$. Note
that $\<\tilde X,v_g\>_1$ is just a multiple of $\<\tilde X,v_g\>$ if
one of the components of $\tilde X = \tilde X_\p + \tilde X_\k$ are zero.
Now we discuss which of the zero curvature planes in 
$G = Sp(2)$ (see Lemma \ref{lemma0}) can be horizontal at any $g\in G$.
By a slight abuse of language, a plane $\tilde\sigma$ spanned
by $\tilde X,\tilde Y \in\g$ will be called  {\it horizontal at $g$} if 
\beq
	\<\tilde X,v_g\>_1 = \<\tilde Y,v_g\>_1 = 0
\eeq
for all $v \in \Im\H$.

\medskip\ms
{\bf Case 1.}

\begin{lemma} \label{lemma1} A plane of type (\ref{XY1})
is nowhere horizontal.
\end{lemma}

\proof 
$\<\tilde Y,v_g\> = \<y,\bar ava-v\> = \<ay\bar a -y,v\>$ vanishes for all
$v\in\Im\H$ iff $y = ay\bar a$, and likewise $\<\tilde X,v_g\>$ vanishes
for all $v$ iff $y = by\bar b$. But this implies $|a|=|b| = 1$ in
contradiction to (\ref{a2+b2}).
\endproof

\medskip\ms
{\bf Case 2a.}

\begin{lemma} \label{lemma2}
If a plane of type (\ref{XY2}) is horizontal at $g$
then either $a=0$ or $b=0$ or
\beq \label{det}
	\det(I-\Ad(a^{-1}) - \Ad(b^{-1})) = 0.
\eeq
\end{lemma}

\proof
The matrix $\tilde X$ is horizontal at $g$ if and only if
\beq \label{bxa}
	0 = \< \tilde X, v_g\> = 2\<x,\bar b v a\> = 2\<bx \bar a,v\>
\eeq
for all $v \in \Im \H$.  This is equivalent to $bx \bar a \in \R$.  
Hence, either $a=0$ or $b=0$ or $b x = r a$ for some non-zero 
$r \in \R$.  In the latter case we have, in particular 
\begin{eqnarray} \label{Adbx}
	\Ad (bx) &=& \Ad (a), \\  
\label{Adx}
	\Ad(x) &=& \Ad(b^{-1})\Ad(a),
\end{eqnarray}
provided that $b\neq 0$.
On the other hand, the matrix $\tilde Y$ is horizontal at $g$ if and only if
\beq \label{vgY}
	0 = \<\tilde Y, v_g\> 
	= \<|a|^2 \Ad (a) y - y + |b|^2 \Ad (bx) y - \Ad (x) y,v\>
\eeq
for all $v \in \Im \H$. Since $y \in \Im\H$, this means
\begin{eqnarray} \label{aby}
	0 &=& |a|^2 \Ad (a) y + |b|^2 \Ad (bx) y - y - \Ad (x) y \\
	  &{\buildrel(\ref{Adbx})\over =} &  \Ad(a)y - y - \Ad(x)y	\cr
	  &{\buildrel(\ref{Adx})\over =}&  \Ad(a)y - y - \Ad(b^{-1})\Ad(a)y
	  \nonumber
\end{eqnarray}
where we have also used $|a|^2+|b|^2 = 1$ (\ref{a2+b2}). If $a\neq 0$,
we obtain from the last equality
$$
	\Ad (a) y \in \ker(I - \Ad (a^{-1}) - \Ad (b^{-1})) 
$$
and in particular 
$$ 
	\det (I - \Ad (a^{-1}) - \Ad (b^{-1}) = 0.
	\eqno{(\ref{det})}
$$
\endproof

\begin{lemma} \label{lemma3} There exists a plane of type (\ref{XY2})
which is horizontal at $g$ if and only if either (\ref{det}) holds or
\beq \label{a=0}
a = 0,\ \ \ |\Im b\,| \geq \2 \ \ \textrm{ or }
\ \ b = 0,\ \ \ |\Im a\,| \geq \2.
\eeq
\end{lemma}

\proof
Suppose first $a,b \neq 0$. If (\ref{det}) is satisfied, there is 
a non-zero $w \in \ker(I - \Ad (a^{-1}) - \Ad (b^{-1}))$. Then defining 
$y = \Ad (a^{-1}) w$ and 
$x = b^{-1} a$, we obtain a horizontal plane of type (\ref{XY2}) at $g$.
The converse conclusion was done before.

\smallskip\noindent
Now suppose $b=0$. Then $|a| = 1$ and  Equation (\ref{aby}) becomes
\beq \label{aday}
	\Ad (a) y - y = \Ad (x) y.
\eeq
Geometrically, this equality means  that $\Ad(a)$ rotates $y$ by the angle 
$\frac{\pi}{3}$ (the three vectors $\,\Ad(a)y,\,y,\,\Ad(x)y\,$ form the sides of an
equilateral triangle). Hence (\ref{aday}) has a solution $(x,y)$ if and only if
the rotation angle of the rotation $\Ad(a)$ is $\geq \frac{\pi}{3}$. This  in
turn is  equivalent to $\sphericalangle(a,1) \geq \frac{\pi}{6}$, i.e. $|\Im a| \geq
\2$.   Inserting the solution $(x,y)$ into (\ref{XY2}) defines a horizontal
plane of type (\ref{XY2}). The case $a=0$ is similar. \endproof

\medskip\goodbreak\noindent
{\bf Case 2b.}

\begin{lemma} \label{lemma4}
If a plane of type (\ref{XY3}) is horizontal at $g$, then
\beq \label{a=b}
	|a| = |b| = 1/\sqrt2
\eeq
and $w := \Im a^{-1}b$ satisfies
\beq \label{w-2awa}
	\<w - 2 a^{-1}wa,w\> = 0.
\eeq
\end{lemma}

\proof
\begin{eqnarray} \label{vgY+}
\<v_g,Y_+\> &=& \<\bar ava+\bar bvb - 2v,y\> 
		= \<v,ay\bar a + by\bar b - 2y\> \\
		\label{vgY-}
\<v_g,Y_-\> &=& \ \ \ \<\bar ava-\bar bvb,y\> \ \ \ \ 
		= \ \ \<v,ay\bar a - by\bar b\> 
\end{eqnarray}
Thus $\<\tilde Y,V_g\> = 0$ iff one of the following equations holds:
\begin{eqnarray} 
	ay\bar a + by\bar b &=& 2y,	\cr
	ay\bar a - by\bar b &=& 0.	\nonumber
\end{eqnarray}
The first of these equations is impossible by the triangle
inequality together with (\ref{a2+b2}):
$$
	|ay\bar a + by\bar b| \leq |ay\bar a| + |by\bar b| 
	\leq (|a|^2+|b|^2)|y| = |y| < |2y|.
$$
Thus we are left with the second equation, 
\beq \label{aybara}
	ay\bar a = by\bar b,
\eeq 
which implies $|a| = |b|$.\\
Note that we have also shown that $Y_+$ cannot be horizontal. Thus we 
need only consider $\tilde X = X_-$ and $\tilde Y = Y_-$ in (\ref{XY3}), 
and 
\beq \label{xy}
	xyx^{-1} = -y
\eeq
which means that $x$ is imaginary and nonzero with $x \perp y$.\\
Now let $\tilde X,\tilde Y$ be as above spanning $\tilde\sigma$. 
By the preceding remark we have
\beq \label{XY}
	\tilde Y = \bpm y & \cr & -y \epm,\ \ \ 
	\tilde X = \bpm y & x \cr x & y \epm
\eeq
with $y \perp x \in \Im\H$. Thus according to (\ref{s})
we get for all $v \in \Im\H$
\begin{eqnarray} \label{vgX}
0 = \<\tilde X,v_g\>_1 &=& 
	2\<x,\bar bva\>\ +\ \tilde s\<y,\bar ava+\bar bvb - 2v\> 
	\cr
	&=& 2\<bx\bar a,v\> + \tilde s\<ay\bar a+by\bar b - 2y,v\>
	\cr
	&=& \<bxa^{-1} + \tilde s(aya^{-1}- 2y),v\>, 	 
\end{eqnarray}
where we have used $2\bar a = a^{-1}$ 
and $ay\bar a = by\bar b = \2 aya^{-1}$ 
from  (\ref{a=b}) and (\ref{aybara}).
Putting $p = a^{-1}b/\tilde s$, we obtain
\beq \label{imbxa-1}
	\Im apxa^{-1} = 2y - aya^{-1}.
\eeq 
From $aya^{-1} = byb^{-1}$ we see $yp = py$, thus 
$p  \in \C_y : = \R + \R y$
and thus the left multiplication with $p$ preserves 
$\C_y$ and $\C_y^\perp$.
By (\ref{xy}) we have $x \in \C_y^\perp$ and therefore 
$px \in \C_y^\perp$. 
Conjugating (\ref{imbxa-1}) by $a^{-1}$ we obtain 
\begin{eqnarray} \label{impxperpy}
	2a^{-1}ya-y &=& \Im(px) \perp y ,	\\
	 \label{y-2aya}
	\<2a^{-1}ya-y,y\> &=& 0.
\end{eqnarray}
Since $w = \Im \tilde sp \in \C_y$ is a multiple of $y$,
we may replace $y$ by $w$ in Equation (\ref{y-2aya})
and obtain (\ref{w-2awa}). 
\endproof

\goodbreak 
\noindent{\bf Remark 1.} 

\hskip 4.7cm
\includegraphics{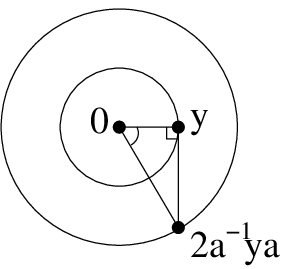}

\ms Geometrically, (\ref{y-2aya}) means that the
angle between $y$ and $a^{-1}ya$ is $\pi/3 = 60^o$. Thus the
rotation angle of $\Ad(a^{-1})$ (and of $\Ad(b^{-1})$,
see (\ref{aybara})) must be $\geq \pi/3$, hence 
$\sphericalangle(1,a) \geq \pi/6$, or in other words, 
\beq \label{ima/a}
	\frac{|\Im a|}{|a|} \geq \2.
\eeq

\begin{lemma} Suppose that $a,b \in \H$ satisfy {\rm (\ref{a=b})}, {\rm
(\ref{w-2awa})} and {\rm (\ref{ima/a})}. Then there exists a horizontal plane
of type (\ref{XY3}) at $g = (\bsm a & b \cr c & d \esm)$.   
\end{lemma}

\proof First suppose that $\tilde p = a^{-1}b = \tilde s p$ is real which in view of
(\ref{a=b})  means $a=\pm b$. By (\ref{ima/a}), the rotation angle of
$\Ad(a^{-1})$ is $\geq \pi/3$, hence there exists a nonzero $y \in
\Im\H$ which is rotated precisely by the angle $\pi/3$ and thus
satisfies (\ref{y-2aya}). Put $x = 2a^{-1}ya-y \perp y$ 
and define $\tilde X,\tilde Y$ 
as in (\ref{XY}). This matrix pair is of type (\ref{XY3}), and it is 
perpendicular to $V_g$  by (\ref{vgY-}) and
(\ref{vgX}).

Now suppose that $w = \Im\tilde p \neq 0$; in this case (\ref{w-2awa})
implies  (\ref{ima/a}). Then we choose $y = w$ and $x =
\Im\left(p^{-1}(2a^{-1}wa-w)\right)$, compare (\ref{impxperpy}). Since
$w-2a^{-1}wa \in \C_y^\perp$ (it is imaginary and perpendicular to $w =
y$),  we also have $p^{-1}(w-2a^{-1}wa) \in \C_y^\perp,$ hence $x\perp
y$ and thus $xyx^{-1} = -y$. Defining matrices $\tilde X,\tilde Y$ using (\ref{XY}), 
these are of type (\ref{XY3}) and
perpendicular to $V_g$ by (\ref{vgY-})  and (\ref{vgX}).
\endproof

\noindent {\bf Remark 2.} Clearly, the relations (\ref{det}), (\ref{a=0}), 
(\ref{a=b}), (\ref{w-2awa}) and
(\ref{ima/a}) must be invariant under the action of $U$. In fact, if
$u = \left((\bsm q & \cr & 1\esm),(\bsm q&\cr&q \esm)\right)$, we have
$u .g = \tilde g = \left(\bsm \tilde a & \tilde b \cr \tilde c & \tilde d
\esm\right)$ with $\tilde a = qaq^{-1}$ and $\tilde b = qbq^{-1}$.

\medskip Now we have proved the following

\begin{theorem} \label{theorem}
Let $G = Sp(2)$ with the left invariant metric {\rm (\ref{t})}
and $U \subset G\x G$ defined by {\rm (\ref{U})}. The orbit space $M = G/U$ inherits a
Riemannian metric such that the canonical projection $\pi : G \to M$ is a
Riemannian submersion. Let 
$$
	Z = \{p \in M;\ \exists{\sigma \subset T_pM}: \ \sec(\sigma) = 0\}.
$$
Then $Z = Z_1\cup Z_2\cup Z_3\cup Z_4$ where 
\goodbreak
\begin{eqnarray}
\pi^{-1}Z_1 &=& \{\left(\bsm a & b \cr c & d \esm\right);\ a,b \neq 0, \
	\det(I-\Ad(a^{-1})-\Ad(b^{-1})) = 0\}, \cr
\pi^{-1}Z_2 &=& \{\left(\bsm a & b \cr c & d \esm\right);\ |a|=|b|,\
	w := \Im a^{-1}b \perp w-2a^{-1}wa, \ 	
	|\Im a| \geq |a|/2\},	\cr
\pi^{-1}Z_3 &=& \{\left(\bsm
a & b \cr c & d \esm\right);\ b = c = 0,\ |\Im a| \geq 1/2\},	\cr
\pi^{-1}Z_4 &=& \{\left(\bsm
a & b \cr c & d \esm\right);\ a = d = 0,\ |\Im b| \geq 1/2\}, \nonumber	
\end{eqnarray}
where all matrices $\left(\bsm a & b \cr c & d \esm\right)$ are supposed 
to belong to $Sp(2)$.
\qed
\end{theorem}

\noindent {\bf Remark 3.} The mistake in \cite{E2} is in the third line of the
proof  of the Theorem,	page 1166. The computation of $\<v_g,X\>$  holds only
for $X\in\k$, but $X$ may have a nonzero $\p$-component as  well. Thus the
matrix $X$ in (4), p. 1166, is too special and must be  replaced with the more
general $X = \left(\bsm  ry & -\bar x \cr x &  -rxyx^{-1} \esm\right)$  for
arbitrary $r\in\R$, and instead of (5)  $\Im(bx\bar a) =0$  we obtain $(5')$
$\Im(bx\bar a) = r(y-ay\bar a)$,   while Equation (6) 
($ay\bar a - y + bxyx^{-1}\bar b - xyx^{-1} = 0$) 
remains unchanged. We have 15
variables, $(a,b) \in S^7$,  $x \in \H$, $y \in \Im(\H)$, $r\in\R$, with two
arbitrary real constants  (the lengths of $x$ and $y$), and 6 constraint 
equations $(5')$ and (6)  which reduce the number of free  variables to 7. Thus
the solution set  is likely to project  onto a subset with positive measure in
the  $(a,b)$-space $S^7$; this would imply that the metric considered in 
\cite{E2} fails to have almost positive curvature.

\bigskip

\end{document}